\documentclass[12pt]{article}

\usepackage{amssymb}
\usepackage{amsthm}
\usepackage{amsmath}
\usepackage[T1]{fontenc}
\usepackage[utf8]{inputenc}
\usepackage[english]{babel}
\usepackage[pdftex]{graphicx}
\usepackage{color}

\usepackage{indentfirst}              
\usepackage{makeidx}                    
\usepackage[nottoc]{tocbibind}         
\usepackage{courier}                   
\usepackage{type1cm}                    
\usepackage{listings}                
\usepackage{titletoc}
\usepackage[all,cmtip]{xy}

\newcommand{\Z}{\mathbb{Z}}

\newtheorem{thm}{Theorem}[section]
\newtheorem{prop}[thm]{Proposition}
\newtheorem{lem}[thm]{Lemma}
\newtheorem{Cor}[thm]{Corollary}

\newtheorem{rem}[thm]{Remark}

\def\bee{\begin{eqnarray}}
\def\bes{\begin{eqnarray*}}
\def\eee{\end{eqnarray}}
\def\ees{\end{eqnarray*}}

\def\e{\varepsilon}
\def\Z{\mathbf Z}

\def\12{\tfrac12}
\def\Proof{{\it Proof.\ }\ }
\def\ctd{\hfill$\Box$}

\def\bee{\begin{eqnarray}}
\def\bes{\begin{eqnarray*}}
\def\eee{\end{eqnarray}}
\def\ees{\end{eqnarray*}}

\def\Reg{\rm{Reg\,}}

\def\e{\varepsilon}

\def\Proof{{\sl Proof.}\ }

\def\0{_{\bar 0}}
\def\1{_{\bar 1}}

\title{A Coordinatization Theorem for the Jordan algebra of symmetric $2\times 2$ matrices}

\author{Jes\'us Laliena   \footnote {Supported by a grant PID 2021-123461 NB-C21 founded by MCIN/AEI/10.13039/5011000110033 and by "ERDF A way of making Europe".} 
\\{\small Dep. de Matem\'aticas y Computaci\'on}\\
{\small  Universidad de La Rioja}\\
{\small  26004, Logro\~no. Spain}\\
{\small jesus.laliena@dmc.unirioja.es }
\and
Victor L\'opez Sol\'{\i}s\\
{\small Dep. Acad\'emico de Ciencias Básicas y Afines }\\
{\small   Universidad Nacional de Barranca}\\
{\small Barranca. Peru}\\
{\small vlopez@unab.edu.pe}
\and
Ivan Shestakov  \footnote {Supported by  FAPESP  grant 2018/23690-6,  by  CNPq  grant 304313/2019-0, and also by IMC SUSTech, Shenzhen,  China.} 
\\{\small Instituto de Matem\'atica e Estat\'{\i}stica}\\
{\small Universidade de S\~ao Paulo}\\
{\small S\~ao Paulo. Brasil}\\
{\small ivan.shestakov@gmail.com}
}

\date{\quad}

\begin{document}
\maketitle\vspace{-1.5cm}

\begin{abstract}
The Jacobson Coordinatization Theorem describes the structure of unitary Jordan algebras containing the algebra $H_n(F)$ of symmetric $n\times n$ matrices over a field $F$ with the same identity element, for $n\geq 3$.
In this paper we extend the Jacobson Coordinatization Theorem for $n=2$. Specifically, we prove that if $J$ is a unitary Jordan algebra containing the Jordan matrix algebra $H_2(F)$  with the same identity element, then $J$ has a form
$J=H_2(F)\otimes A_0+k\otimes A_1$, where $A=A_0+A_1$ is a $\Z_2$-graded Jordan algebra with a partial odd Leibniz bracket $\{,\}$ and $k=e_{12}-e_{21}\in M_2(F)$, with the multiplication given by
$$
(a\otimes b)(c\otimes d)=ac\otimes bd+[a,c]\otimes \{b,d\},
$$
the commutator $[a,c]$ is taken in $M_2(F)$. 
\end{abstract}

{\parindent= 4em \small  \sl Keywords: Jordan algebra, Jordan matrix algebra, coordinatization theorem.}

\section{Introduction}

In what follows, all algebras are defined over a field $F$ which in case of Jordan algebras is assumed to be  of characteristic different from $2$.

It was proved by H. M. Wedderbun that if  $B$ is an associative algebra with a unitary element $1$, and $A$ is a finite dimensional central simple subalgebra containing $1$, then $B$ is isomorphic to the Kronecker product $A\otimes S$, where $S$ is the subalgebra of elements of $B$ which commute with every element of $A$ .

In 1951 (see \cite {K}), I. Kaplansky proved a similar result for a unitary alternative algebra $B$ and a subalgebra $A$  of $B$ containing $1$ and having the structure of a split Cayley-Dickson algebra. Then $B= A\otimes C$, where $C$ is the center of  $B$ (that is the set of elements of $B$ which commute with every element of $A$ and associate with every pair of elements). N. Jacobson, in 1954 (see \cite{J1}), also proved  a Kronecker factorization theorem for the case when $B$ is a Jordan algebra with $1$ and $A$ is an exceptional simple $27$-dimensional Jordan Albert algebra, which is a subalgebra of $B$ having the same identity as $B$. Then $B= A\otimes C$, where $C$ is the center of $B$. 

In the case of superalgebras, M. L\'opez-D\'{\i}az and I. Shestakov studied in 2002 and 2005 (see \cite{LDS2, LDS1}) the representations of simple alternative and exceptional Jordan superalgebras in characteristic $3$ and, through these representations, they obtained some analogues of the Kronecker factorization theorem for these superalgebras. 
In 2003 (see \cite{CE}), C. Martinez and E. Zelmanov  obtained a Kronecker Factorization Theorem  for the exceptional 10-dimensional Kac superalgebra $K_{10}$.  In 2020 (\cite{PS}), A.Pozhidaev and I.Shestakov proved some analogues of the Kronecker factorization theorem for noncommutative Jordan  superalgebras.

We recall now that, in $1954$, Nathan Jacobson proved the Strong Coordinatization Theorem (see Chapter III in \cite{J}), that is, if $J$ is a unitary Jordan algebra containing a subalgebra $J_0$ isomorphic to the Jordan matrix algebra $H_n(F)$ with $n\geq 3$, and both with the same identity element, then $J\cong H_n(D,*)$, where $D$ is an associative algebra with involution $*$ if $n\geq 4$ and $D$ is an alternative algebra with involution with $*$-symmetric elements  in the nucleus if $n=3$.

For the case $n=2$, M.Osborn \cite{O} classified the simple Jordan algebras with two orthogonal idempotents $e_1,e_2$ whose sum is 1,  under the condition that the Pierce components $J_{11}$ and $J_{22}$ are Jordan division algebras.   

\smallskip

Here we are interested in describing the unital Jordan algebras $J$ such that $J$ contains the Jordan matrix algebra  $H_2(F)$ with the same identity element. 
We do not assume additional conditions on simplicity of $J$ and on the structure of the Pierce components $J_{ii}$.
The analogous problem for unitary alternative algebras containing a split quaternion algebra  with the same identity element has been solved by V. L\'opez-Sol\'{\i}s and I. Shestakov in  $\cite{SS}$ (see also \cite {GS}). 

\smallskip
Let  $\{e_{ij}\,|\, i,j=1,2\}$ be the canonical basis of the $2\times 2$ matrix algebra $M_2(F)$.  We will denote  $e=e_{11}, \,h=e_{12}+e_{21}, \,k= e_{12}-e_{21}$.   Observe that $M_2(F)=H_2(F)+Fk$.  Below we will  call a unital Jordan algebra $J$ containing $H_2(F)$ as a unital subalgebra,  an {\it $H_2$-algebra.} 

Our main result gives  a characterization of  $H_2$-algebras:

{\it Let $J$ be an  $H_2$-algebra over a field  $F$ containing an element $\epsilon$ with $\epsilon ^2=-1$. Then there exists a $\Z_2$-graded unital Jordan algebra $S_0+S_1$, with  a partial anticommutative odd bracket $\{S_i,S_j\}\subseteq S_{i+j+1},\ (i,j)\neq (1,1),$ that satisfies the identities 
\begin{align*}
(*) \quad \quad (x,y,z)&=4\{y,\{x,z\}\}, \text{ for any $x, y, z\in S$,}\\
 \{xz,y\}&=x\{z,y\}+z\{x,y\}=\{z,xy\}+\{x,zy\}, \text{ for any $x, y, z \in S$,}\\
  (\{x,y\},&z,t)+(\{t,x\},z,y)+(\{y,t\},z,x)=0, \text{ for any $x, y\in S_0$ and $z, t\in S_1$},
\end{align*}
 whenever the involved products are defined, such that $J$ is isomorphic to $H_2\otimes S_0+Fk\otimes S_1$, with the multiplication 
\bes
(**) \quad \quad (a\otimes x)(b\otimes y)=a\cdot b\otimes xy+[a,b]\otimes \{x,y\},
\ees
where for $a,b\in M_2(F)$ we denote   $a\cdot b=\tfrac12(ab+ba), [a,b]=ab-ba$.  
Conversely,  for any $\Z_2$-graded unital Jordan algebra $S$ with a partial anticommutative bracket $\{.,.\}$ as above,  the algebra 	$J=H_2\otimes S_0+Fk\otimes S_1$ with the multiplication (**) is a unitary Jordan algebra containing $H_2(F)$ as a unital subalgebra.}

Furthermore,  let $J$ be an $H_2$-algebra,  $J=H_2\otimes S_0+Fk\otimes S_1$,  where $S=S_0+S_1$ is a $\Z_2$-graded Jordan algebra with a partial odd anticommutative bracket defined above.
Then
\begin{itemize}
\item 
{\it The bracket $\{,\}$ is completely defined on $S$ if and only if $J$ is isomorphic to the Jordan algebra of hermitian $2\times 2$ matrices $H_2(A,*)$ over a certain associative algebra with involution $(A,*)$.}
\item
{\it The bracket $\{,\}$ is trivial $(=0)$ if and only if $J$ is isomorphic to a Jordan algebra of bilinear form over an associative commutative algebra $A$.}
\end{itemize}

It remains an open question whether there exist examples of $H_2$-algebras different of hermitian $2\times 2$ matrices or algebras of bilinear form.

\section {The subalgebra $S$.}

A \textit{Jordan algebra} is a vector space $J$ with a bilinear binary operation $(x,y)\mapsto xy $ satisfying the following identities:
\begin{equation}\label{e1}
\begin{split}
xy& = y x,\\
(x^{ 2} y) x& = x^{2} (y x).
\end{split}
\end{equation}
for all $x,y\in J$.  

The following identities holds in every Jordan algebra (see for example \cite{J}, Chapter 1, Section 7)
\bee
(x,y,z)&=&-(z,y,x),\label{id_2}\\
(zt,x,y)+(yz,x,t)+(yt,x,z)&=&0,\label{id_3}\\
(x,yz,t)&=&y(x,z,t)+(x,y,z)t,\label{id_4}
\eee
where $(x,y,z)= (xy)z-x(yz)$.

We remind also the Teichm\"uller identity, which is valid in any algebra:
\bee\label{Teich}
(xy,z,t)+(x,y,zt)=x(y,z,t)+(x,yz,t)+(x,y,z)t.
\eee

We will use the following function introduced by  S.V. Pchelintsev and I.P. Shestakov in \cite{P-S} in Jordan algebras
$$
k(x,y;z,t):=(xy,z,t) -(x,z,t)y-x(y,z,t).
$$
They prove in  \cite{P-S} the following

\begin{lem}\label{lem_1}\cite[Lemma 7.4]{P-S}
The function $k(x,y;z,t)$ is symmetric in $x,y$ and in $z,t$. Moreover, $k(x,y;z,t)=k(z,t;x,y)$.
\end{lem}
\Proof We will give the proof for completeness. It is clear that the function $k$ is symmetric in $x,y$. 
Furthermore,
\bes
k(x,y;z,t)&=&(xy,z,t)-x(y,z,t)-(x,z,t)y\\
&\stackrel{\eqref{Teich}}=&-(x,y,zt)+(x,yz,t)+(x,y,z)t-y(x,z,t)\\
&\stackrel{\eqref{id_2},\eqref{id_3}}=&(zt,y,x)+y(x,z,t)+(x,y,t)z-t(z,y,x)-y(x,z,t)\\
&=&(zt,y,x)-z(t,y,x)-(z,y,x)t=k(z,t;y,x).
\ees
So $k(x,y;z,t)=k(z,t;y,x)=k(t,z;y,x)=k(x,y;t,z).$
\ctd

Now we suppose that $J$ is a unital Jordan algebra and $A$ is a unital subalgebra of $J$ such that $A\cong H_2(F)$. Let $Z=\{z\in J\,|\, (z,a,b)=0$ for all $a,b\in A\}$. 
\begin{lem}\label{lem_2}
$Z$ is a subalgebra of $J$ such that $(Z,J,A)=0$. Moreover, $ZA\cong Z\otimes A$ as vector space.
\end{lem}

\Proof 
 Let $J=J_1+J_{1/2}+J_0$ be the Pierce decomposition of $J$ with respect to $e=e_{11}$ (see for instance \cite{J}, Chapter 3, Section 1).
We prove first that $Z$ is a subalgebra of $J$. It suffices to prove that $(z^2,a,b)=0$ for any $z\in Z,\ a,b\in A$. 
Let $z=z_1+z_{1/2}+z_0$, with $z_0\in J_0, z_1\in J_1, z_{1/2}\in J_{1/2}$. Then the equality
 $(z,e,e)=0$ implies that $z_{1/2}=0$ and  $Z\subseteq J_1+J_0$.
 If $a\in J_0,\, x\in J$, then, because of \eqref{id_2} and \eqref{id_3}, $(a,x,e)=(a,x,e^2)=2(ae,x,e)=0$. Also if $a\in J_1$, then, again because of  \eqref{id_2} and \eqref{id_3}, $(a,x,e)=(a,x,e^2)=2(ae,x,e)=2(a,x,e)=0$, and therefore $(J_0+J_1,J,e)=0$.
 Since $J_1$ and $J_0$ are orthogonal subalgebras, this implies that  $(z^2,J,e)=0$. Observe that we have also proved  that  $(Z,J,e)=0$.
 Futhermore, from $(z,h,h)=0$ we have $z=(zh)h$, therefore, because of  \eqref{id_2} and \eqref{id_3}$,(z,x,h)=((zh)h,x,h)=1/2(zh,x,h^2)=0$ for any $x\in J$,
hence $(z,J,h)=0$ and $(Z,J,A)=0$.   In particular, we have $(z,z,h)=0$. Finally, from \eqref{Teich} and because of $(Z,J,A)=0$, $(z^2,h,h)=-(z,z,h^2)+z(z,h,h)+(z,zh,h)+(z,z,h)h=0$. And the same $(z^2,e,h)=-(z,z,eh)+z(z,e,h)+(z,ze,h)+(z,z,e)h=0$. This proves that $ZZ\subseteq Z$.

Assume now that $z+ue+vh=0$ for some $z,u,v\in Z$. We have $z+ue\in J_1+J_0,\ vh\in J_{1/2}$, hence $z+ue=0,\ vh=0$. Therefore, $v=(vh)h=0$. Furthermore, $ue=-z\in Z$, hence $0=(ue,e,h)=(ue) h-\tfrac12 (ue)h=\tfrac12 uh -\tfrac14 uh=\tfrac14 uh$. Now $0=(uh)h=u$ and $z=0$.

\ctd  

Consider $J$ as $A$-bimodule, then by \cite{J} (Sections 11 and 13 in Chapter 2, and also by \cite{J2}), $J=M\oplus N$ as a vector space, where $M$ is a direct sum of copies of the regular bimodule  $\Reg H_2(F)$ and $N$ is a direct sum of one-dimensional bimodules where $e|_N=\tfrac12,\ h\cdot N=0$.

\begin{lem}\label{lem_3}
We have $M=AZ,\ N^2\subseteq Z,\ ZN\subseteq N$.
\end{lem}
\Proof For any $z\in Z$ we have $zA\cong \Reg A$, hence $ZA\subseteq M$. On the other hand, if an $A$-submodule $V\subseteq J$ is isomorphic to $\Reg A$, then $V=vA$ where $v\in Z$, hence $V\subseteq ZA$.
Furthermore, for any $n\in N,\ x\in J$ we have, because of \eqref{id_2} and \eqref{id_3},
\bes
(n^2,x,h)&=&2(n,x,nh)=0,\\
(n^2,x,e)&=&2(n,x,ne)=(n,x,n)=0.
\ees
Therefore, $(n^2,J,A)=0$ and $n^2\in Z$. By commutativity, $N^2\subseteq Z$.

Finally, for any $n\in N,\,z\in Z$ we have 
\bes
zn&\in& (J_1+J_0)J_{1/2}\subseteq J_{1/2},\\
(zn)h&=&(z,n,h)=(zh\cdot h,n,h)=\tfrac12(zh,n,h^2)=0.
\ees
Therefore, $V=F(zn)$ is a one-dimensional $A$-submodule such that $R_e|_V=\tfrac12,\ R_h|_V=0$, hence $V\subseteq N$ and $ZN\subseteq N$.

\ctd

\begin{Cor}\label{cor_1}
$S=Z+N$ is a $\Z_2$-graded subalgebra of $J$ with $S_0=Z,\ S_1=N$.
\end{Cor}

Let now try to calculate other products in $J$.

\begin{prop}\label{prop_1}
Let $z,v\in Z,\, n\in N$. then
\bee
(ez)(ev)&=&e(zv),\label{id_6}\\
(hz)(hv)&=&zv,\label{id_7}\\
(ez)(hv)&=&\tfrac12 h(zv)+\{z,v\},\label{id_8}\\
(ez)n&=&\tfrac12 zn+[z,n]h,\label{id_9}\\
(hz)n&=&2[z,n](1-2e),\label{id_10}
\eee
where $\{z,v\}=(ev,z,h)\in N,\ [z,n]=-[n,z]=h(z,e,n)\in Z$ are bilinear functions that satisfy the identities
\bee
\{z,v\}&=&-\{v,z\},\  \{z,v^2\}=2v\{z,v\}=2\{zv,v\},\label{id_11}\\ \ 
[z^2,n]&=&2z[z,n]=2[z,zn],\label{id_12}\\
\{z,n^2\}&=&-2n[z,n].\label{id_13}
\eee
\end{prop}
\Proof
Since $(Z,J,A)=0$ from Lemma \ref{lem_2}, for any $a,b\in A$ we have
\bes
(az)(bv)=(az\cdot v)b -(az,v,b)=(a\cdot zv)b-(az,v,b)=(ab)(zv)-(az,v,b).
\ees
Furthermore,
\bes
(az,v,a)=-\tfrac12 (a^2,v,z)=0,
\ees
hence $(az,v,b)$ is antisymmetric on $a,b$ and the only non-zero case for $a,b\in \{1,e,h\}$ is $(hz,v,e)=-(ez,v,h)$, because of \eqref{id_3}.  We will denote this as $\{z,v\}$.
Furthermore,  any commutative algebra satisfies the identity
\bee\label{id_cicl}
 (x,y,z)+(y,z,x)+(z,x,y)=0,
 \eee 
therefore, because of \eqref{id_cicl}, \eqref{id_2} and \eqref{id_4},
\bes
(ez,v,h)=-(v,h,ez)-(h,ez,v)=(ez,h,v)-e(h,z,v)-z(h,e,v)=(ez,h,v).
\ees
Since $e,z,v\in J_0+J_1$ and $h\in J_{1/2}$, it is clear that $(ev,z,h)=(ev,h,z)\in J_{1/2}$. Moreover, from \eqref{id_4},  $(ev,h,z)h=\tfrac12 (ev,h^2,z)=0$, hence $\{z,v\}\in N$.
And also, from \eqref{id_3}, it follows that 
\bes
\{z,v\}= -(ez,v,h)=-(ez,h,v)=(zv,h,e)+(ev,h,z)=(ev,z,h).
\ees
Now, using \eqref{id_cicl}, \eqref{id_2}, \eqref{id_3}, \eqref{id_4} and that $(Z,J,A)=0$, we have
\bes
\{z,z\}&=&(ez,z,h)= -(z,h,ez)-(h,ez,z)=(ez,h,z)-e(h,z,z)-z(h,e,z)=0,\\
\{z,v^2\}&=&(hz,v^2,e)=2v(hz,v,e)=2v\{z,v\},\\
\{z,v^2\}&=&-(ez,v^2,h)=-(ez,h,v^2)=-2((ez)v,h,v)=-2(e(zv),h,v)=2\{zv,v\}.
\ees

Furthermore, because of \eqref{id_2} and \eqref{id_3}, 
\bes
(e,x,[z,n])&=&(e,x,h(z,e,n))=((z,e,n)e,x,h)+(eh,x,(z,e,n))\\
&=&\tfrac12((z,e,n),x,h)+\tfrac12(h,x,(z,e,n))=0,\\
(h,x,[z,n])&=&(h,x,h(z,e,n))=\tfrac12(h^2,x,(z,e,n))=0.
\ees
thus $(A,J,[z,n])=0$ and $[z,n]\in Z$.

Now, using again \eqref{id_2} and \eqref{id_3}, 
\bes
[z^2,n]&=&h(z^2,e,n)=2h(z,e,zn)=2[z,zn].
\ees 
From \eqref{id_2} and Lemma \ref{lem_1} and \eqref{id_cicl} we observe that 
\bes
[z,n]&=&h(z,e,n)=-h(n,e,z)=-(hn,e,z)+n(h,e,z)+k(h,n;e,z)\\
&=&k(h,n;z,e)=-h(n,z,e)=h(e,z,n),
\ees
therefore, because of \eqref{id_4},
\bes
[z^2,n]&=&h(e,z^2,n)=2h(z(e,z,n))=2z(h(e,z,n))+2(z,(e,z,n),h)\\
&=&2z[z,n]+2(z,(e,z,n),h).
\ees
But $(z,J,h)=0$, hence  $[z^2,n]=2z[z,n]$.
\smallskip

Now, $(ez)n\in J_{1/2}$, hence $(ez)n=z_1h+n_1$ for some $z_1\in Z,\, n_1\in N$. Multiplying by $h$ this equality, we have
$(ez\cdot n)h=(z_1\cdot h)h=z_1$, and by \eqref{id_3} and Lemma \ref{lem_3},
\bes
n_1&=&ez\cdot n-((ez\cdot n)h)h=-(ez\cdot n,h,h)=(ez\cdot h,h,n)\\
&=&(z\cdot eh,h,n)=\tfrac12 (zh,h,n)=\tfrac12 (h,h,zn)=\tfrac12 zn.
\ees
Again because of Lemma \ref{lem_3} we observe that $z_1=(ez\cdot n)h=((e,z,n)+e(zn))h=(e,z,n)h+\tfrac12 (zn)h$, and that  $(zn)h\in Nh=0$, hence, from \eqref{id_cicl} and \eqref{id_2}, $z_1=(e,z,n)h=(z,e,n)h=[z,n]$.
Therefore, $(ez)n=\tfrac12 zn+[z,n]h$.

\smallskip

Furthermore, $(hz)n\in J_1+J_0$, hence $(hz)n=(z,h,n)=z_1+z_2e$ for some $z_1,z_2\in Z$. Since, because of \eqref{id_4},  $(z,h,n)h=\tfrac12 (z,h^2,n)=0$, we have $z_1h+\tfrac12 z_2h=0$, which gives $z_2=-2z_1$ and $(z,h,n)=z_1(1-2e)$. Since $(1-2e)^2=1$, we have, because of \eqref{id_4}, $z_1=(z,h,n)(1-2e)=(z,h,n)-2(z,eh,n)+2h(z,e,n)=2h(z,e,n)=2[z,n]$.

\smallskip

Finally, consider $\{n^2,z\}=-\{z,n^2\}=-(hz,n^2,e)=-2n(hz,n,e)$, because of \eqref{id_4}. We have $(hz,n,e)=(hz\cdot n)e-\tfrac12(hz)n=(hz\cdot n)(e-\tfrac12)=-\tfrac12 z_1=-[z,n]$, hence 
$\{n^2,z\}=-2n[n,z]$, proving \eqref{id_13}.
\ctd

\begin{prop}\label{prop_2}
The operations $\{,\}$ and $[,]$ satisfy the identities
\bee
4[w,\{u,v\}]&=&(u,w,v),\label{id_15} \\
4\{w,[u,n]\}&=&(u,w,n),\label{id_16} \\
4[[u,m],n]&=&(u,n,m),\label{id_17} \\ 
\ [w,un]&=&[uw,n]+[nw,u] \ \label{id_18} \\
&=&u[w,n]-n\{w,u\},\label{id_19} \\ 
\ [k,mn]&=&[km,n]+[kn,m],\label{id_20}
\eee
for any $u,v,w\in Z,\ k,m,n\in N$.
\end{prop}
\Proof
Consider, for $u,v,w\in Z$, because of \eqref{id_3}, \eqref{id_8}, and because of Lemma \ref{lem_2} and Proposition \ref{prop_1},
\bes
0&=&(eu\cdot w,h,hv) +(eu\cdot hv,h,w)+(hv\cdot w,h,eu)\\
&=&(e\cdot uw,h,hv)+(\tfrac12 h\cdot uv+\{u,v\},h,w)+(h\cdot vw,h,eu)\\
&=&\tfrac12 (h\cdot uw)(hv)-(e\cdot uw)v+\tfrac12(uv) w-\tfrac12 (uv)w\\
&-&2[w,\{u,v\}](1-2e)+(vw\cdot u)e-\tfrac12 (vw)u\\
&=&\tfrac12 (uw)v-e(uw\cdot v)-2[w,\{u,v\}](1-2e)+(vw\cdot u)e-\tfrac12 (vw)u\\
&=& -2[w,\{u,v\}](1-2e)+u,w,v)(1/2-e).
\ees
Since $[w,\{u,v\}], (u,w,v)\in Z$,  multiplying by $1-2e$ and using Lemma \ref{lem_2} we obtain \eqref{id_15}.

\smallskip

Furthermore, considering, and using \eqref{id_3}, Proposition \ref{prop_1}, Lemma \ref{lem_2} and Lemma \ref{lem_3}, 
\bes
0&=&(we\cdot uh,h,n)+(we\cdot n,h,uh)+(uh\cdot n,h,we)\\
&=& (\tfrac12 h\cdot wu+\{w,u\},h,n)+(\tfrac12wn+[w,n]h,h,uh)+(2[u,n](1-2e),h, we)\\
&=& \tfrac12 (wu)n-\tfrac12(wn)u+[w,n](uh)-[w,n]h\cdot u-2[u,n](1-2e)\cdot \tfrac12(wh)\\
&=&\tfrac12(u,w,n)+(u,h,[w,n])-([u,n]w)h+([u,n]w)h+2\{[u,n],w\},
\ees
we get 
$$
(u,w,n)=4\{w,[u,n]\},
$$ 
proving \eqref{id_16}.

\smallskip
Similarly, again from \eqref{id_3}, Proposition \ref{prop_1} and Lemma \ref{lem_3},
\bes
0&=&(m\cdot uh,n,h)+(h\cdot uh,n,m)+(mh,n,uh)\\
&=&( 2[u,m](1-2e))n\cdot h+(u,n,m)\\
&=&-4(([u,m]e)n)h+(u,n,m)= -4(\tfrac12[u,m]n+[[u,m],n]h)h+(u,n,m)\\
&=&-4[[u,m],n]+(u,n,m)
\ees
implies \eqref{id_17}.
\smallskip

Now, consider
\bes 
[w,un]=h(w,e,un)=h(-(u,e,wn)+(uw,e,n))=-[u,wn]+[uw,n]=[nw,u]+[uw,n],
\ees
because of \eqref{id_3} and Proposition \ref{prop_1}. And this proves  \eqref{id_18}. Similarly,
\bes
[k,mn]=h(k,e,mn)=h((kn,e,m)+(km,e,n))=[kn,m]+[km,n],
\ees
proving \eqref{id_20}. 
For \eqref{id_19}, by \eqref{id_2}, Lemma \ref{lem_1} and Lemma \ref{lem_2}, we have
\bes
[w,un]&=&h(w,e,un)=h(u(w,e,n)+n(w,e,u)-k(u,n;e,w))\\
&=&u(h(w,e,n))+(u,(w,e,n),h)-k(e,w;n,u)h).
\ees
Now using \eqref{id_4}, Lemma \ref{lem_2} and Lemma \ref{lem_3}, we have
\bes
-k(e,w;n,u)h&=&-(ew,n,u)h+(e(w,n,u)+w(e,n,u))h\\
&=&-(ew,nh,u)+(ew,h,u)n+\tfrac12(w,n,u)h\\
&=&(ew,h,u)n+n'h=(ew,h,u)n\\
&=&\tfrac12(wu\cdot h)n-\tfrac12(wu\cdot h)n-\{w,u\}n\\
&=&-\{w,u\}n,
\ees
with $n^{\prime} \in N$.

Finally,
\bes
(u,(n,e,w),h)=(u, \tfrac12 nw-\tfrac12nw-[w,n]h,h)=-(u,[w,n]h,h)=0,
\ees
proving  \eqref{id_19}.

\ctd

Denote by $\{,\}$ a partial bracket on $S=Z+N$, defined by
\bes
\{x,y\}=\left\{\begin{array}{c} 
\{x,y\},\ x,y\in Z,\\ \ 
\ \ -[x,y],\  \ x\in Z,\, y\in N.
\end{array}\right.
\ees
Observe that partial bracket $\{,\}$ is an odd bracket with respect to the $\Z_2$-grading of $S$ with $S_0=Z,\, S_1=N$; that is, $\{S_i,S_j\}\subset S_{i+j+1\hbox{\tiny{ (mod 2) }}}$. 
\begin{Cor}\label{cor_2}
The $\Z_2$-graded Jordan algebra $S=Z+N$ with the partial odd bracket $\{,\}$ satisfies the identities
\bee
(x,y,z)&=&4\{\{x,z\},y\},\label{id_21}\\
\{xy,z\}&=&x\{y,z\}+y\{x,z\}=\{y,xz\}+\{x,yz\},\label{id_22}
\eee
always when all the involved products are defined.
\end{Cor}
\Proof
It follows from identities \eqref{id_11} - \eqref{id_13}, \eqref{id_15} - \eqref{id_20}, and their linearizations.

\ctd

We will need also
\begin{lem}\label{lem_2.33}
The following identity is verified in $S=Z+N$ 
\bee\label{id_22'}
(\{z,u\},m,n)+(\{u,n\},m,z)+(\{n,z\},m,u)=0,
\eee
for any $z,u\in Z,\ m,n\in N$.
\end{lem}
\Proof
For any $z,u\in Z,\ n,m\in N$ we have due to \eqref{id_3}
\[
S=(ez\cdot hu,m,n)+(ez\cdot n,m,hu)+(hu\cdot n,m,ez)=0.
\]
Consider every summand of this sum, and from Proposition \ref{prop_1} and Lemma \ref{lem_3}
\bes
(ez\cdot hu,m,n)&=&(\tfrac12 h(zu)+\{z,u\},m,n)\\
&=&([zu,m](1-2e)+\{z,u\}m)n-(\tfrac12h\cdot zu+\{z,u\})\cdot mn\\
&=&-2[[zu,m],n]h+(\{z,u\},m,n)-\tfrac12h(zu\cdot mn),
\ees
\bes
(ez\cdot n,m,hu)&=&(\tfrac12 zn+[z,n]h,m,hu)\\
&=&(\tfrac12 zn\cdot m+2[[z,n],m](1-2e))\cdot hu\\
&-&(\tfrac12 zn+[z,n]h)(2[u,m](1-2e))\\
&=&\tfrac12((zn\cdot m)u)h-4\{[[z,n],m],u\}\\
&+&2[[u,m],zn]h+4\{[u,m],[z,n]\},
\ees
\bes 
(hu\cdot n,m,ez)&=& (2[u,n](1-2e),m,ez)\\
&=&-4[[u,n],m]h\cdot ez-(2[u,n]-4[u,n]e)(\tfrac12 zm+[z,m]h)\\
&=&2([[u,n],zm]-[[u,n],m]z)h\\
&+&4(\{[u,n],[z,m]\}-\{z,[[u,n],m]\}).
\ees
Therefore,
\bes
&S=&(-2[[zu,m],n]-\tfrac12(zu\cdot mn)+\tfrac12((zn\cdot m)u)+2[[u,m],zn]+2([[u,n],zm]-[[u,n],m]z))h\\
&+&(\{z,u\},m,n)-4(\{[[z,n],m],u\}-\{[u,m],[z,n]\}-\{[u,n],[z,m]\}+\{z,[[u,n],m]\})\\
&=&S_1h+(\{z,u\},m,n)-S_2.
\ees
We have
\bes
S_1&=&\tfrac12((zn\cdot m)u-zu\cdot mn)+2(-[[zu,m],n]+[[u,m],zn]+[[u,n],zm]-[[u,n],m]z)\\
&=&\tfrac12((mn,z,u)-(m,n,z)u)+2(-[[zu,m],n]+[[u,m],zn]+[[u,n],z]m)\\
&\stackrel{\eqref{id_21}}=&\tfrac12((mn,z,u)-(m,n,z)u -(zu,n,m)+(u,zn,m)  +(u,z,n)m )\\
&=&\tfrac12(-(u,z,nm)+u(z,n,m)-(uz,n,m)+(u,zn,m) +(u,z,n)m )\stackrel{\eqref{Teich}}=0.
\ees
\bes 
S_2&=&4(\{[[z,n],m],u\}-\{[u,m],[z,n]\}-\{[u,n],[z,m]\}+\{z,[[u,n],m]\})\\
&\stackrel{\eqref{id_16}}=&-([z,n],u,m)+(u,[z,n],m)-(z,[u,n],m)+([u,n],z,m)\\
&\stackrel{\eqref{id_cicl}}=&([u,n],m,z)-([z,n],m,u).
\ees
Hence we have
\bes
0&=&(\{z,u\},m,n)-([u,n],m,z)+([z,n],m,u)\nonumber\\
&=&(\{z,u\},m,n)+(\{u,n\},m,z)+(\{n,z\},m,u),
\ees
proving the lemma.

\ctd

Observe that if the bracket $\{.,.\}$ is defined on the whole algebra $S$ then identity \eqref{id_22'} follows from identities \eqref{id_21},\eqref{id_22}. In fact, in this case we have
for any $x,y,z,t\in S$
\bes
&(\{x,y\},z,t)+(\{t,x\},z,y)+(\{y,t\},z,x)=&\\
&\stackrel{\eqref{id_21}}=4\{ \{\{x,y\},t\}+\{\{y,t\},x\}+\{\{t,x\},y\},z\}&\\
&=\{(x,t,y)+(t,y,x)+(y,x,t),z\}\stackrel{\eqref{id_cicl}}=0.&
\ees

We have proved
\begin{thm}\label{thm_1}
Let $J$ be a unital Jordan algebra containing $A=H_2(F)$ as a unital subalgebra. Then there exists a $\Z_2$-graded Jordan algebra $S=Z+N, \ (S_0=Z,\,S_1=N),$ with  a partial anticommutative odd bracket 
\bes
\{.,.\}:Z\times Z\rightarrow N,\ \ \ \ \{.,.\}:Z\times N\rightarrow Z,
\ees
that satisfies identities \eqref{id_21} - \eqref{id_22'}, such that $J$ is isomorphic to $A\otimes Z+N$, with the multiplication defined according to the formulas  \eqref{id_6} -- \eqref{id_10} of Proposition \ref{prop_1}, taking in account that $AZ\cong A\otimes Z$ and $[z,n]=-\{z,n\}$.
\end{thm}
\smallskip

Let now $J=H_2(F)\otimes Z+N$ be an algebra constructed from a $\Z_2$-graded Jordan algebra $S=Z+N$ with a partial odd bracket $\{,\}$ satisfying identities \eqref{id_21} - \eqref{id_22} according the construction above. Our next objective will be to prove that the algebra $J$  is Jordan, that is, the Theorem \ref{thm_1} in fact gives a characterization of unital Jordan algebras containing $H_2(F)$ with the same unit.  

\smallskip

Represent the matrix algebra $M_2(F)$ in the form $M_2(F)=H_2(F)+Fk$, where $k=e_{12}-e_{21}$. Then the adjoint Jordan algebra $(M_2(F))^{(+)}$ is a $\Z_2$-graded algebra with $(M_2(F))^{(+)}_{0}=H_2(F)$ and $(M_2(F))^{(+)}_1=Fk$.

\begin{prop}\label{prop_3}
Assume that the field $F$ contains an element $\varepsilon$ with $\varepsilon^2=-1$. The algebra $J=H_2(F)\otimes Z+N$ is isomorphic to the algebra $J'=H_2(F)\otimes Z+Fk\otimes N$, with the multiplication
\bee\label{id_mult}
(a\otimes x)(b\otimes y)=a\cdot b\otimes xy+ \varepsilon [a,b]\otimes \{x,y\},
\eee
where $a,b\in H_2(F)\cup Fk,\, x,y\in Z\cup N;\ a\cdot b,\, xy$ are  Jordan product in $M_2(F)^{(+)}$ and in $S=Z+N$, $[a,b]$ is the commutator in $M_2(F)$, and $\{x,y\}$ is the partial odd bracket in $S$.  
\end{prop} 
\Proof
Observe first that the multiplication is well defined since for $m,n\in N$ we have
\bes
(k\otimes m)(k\otimes n)=k^2\otimes  mn=-1\otimes mn.
\ees 
Consider the application $f:J=A\otimes Z+N\rightarrow J',\ a\otimes z+ n\mapsto a\otimes z+ \varepsilon k\otimes n$, where $a\in A=H_2(F),\, z\in Z,\, n\in N$. It is clear that $f$  is an isomorphism of vector spaces. Let us prove that $f$ is an algebra isomorphism.  We have for $a\in A,\,v,z\in Z,\,n,m\in N$:
\bes
f(a\otimes z\cdot a\otimes v)&=&f(a^2\otimes zv)=(a\otimes z)(a\otimes v)=f(a\otimes z)f(a\otimes v);\\
f(e\otimes z\cdot h\otimes v)&=&f( \tfrac12 h\otimes zv+\{z,v\})=\tfrac12 h\otimes zv+\e k\otimes \{z,v\}\\
&=&eh\otimes zv+\e[e,h]\otimes \{z,v\}=(e\otimes z)(h\otimes v)=f(e\otimes z)f(h\otimes v);\\
f(e\otimes z\cdot n)&=&f(\tfrac12 zn+h\otimes [z,n])=\tfrac12   \varepsilon k\otimes zn+h\otimes [z,n]\\
&=&e\cdot \varepsilon k\otimes zn+\e [e,\e k]\otimes \{z,n\}=(e\otimes z)(\e k\otimes n)=f(e\otimes z)f(n);\\
f(h\otimes z\cdot n)&=&f((1-2e)\otimes 2[z,n])=2(1-2e)\otimes [z,n]=-[h,k]\otimes \{z,n\}\\
&=&(h\otimes z)(\e k\otimes n)=f(h\otimes z)f(n);\\
f(nm)&=& 1\otimes nm=(\e k\otimes n)(\e k\otimes m)=f(n)f(m).
\ees
This proves the proposition.

\ctd

Therefore, it suffices to prove that the algebra $J'=H_2(F)\otimes Z+Fk\otimes N$ is Jordan. 
Let us change the bracket $\{x,y\}$ on $S$ to $\e\{x,y\}$, maintaining the same notation. Then this new bracket still satisfies  identities \eqref{id_22}, \eqref{id_22'}, and instead of \eqref{id_21} it satisfies
\bee\label{id_21'}
(x,y,z)=4\{y,\{x,z\}\}.
\eee
 
Observe also that the commutator $[a,b]$ in $M_2(F)$ is  an odd anticommutative  bracket on the $\Z_2$-graded Jordan algebra $(M_2(F))^{(+)}=H_2(F)+Fk$ that satisfies the identities
\bee
4(a,b,c)_+&=&[b,[a,c]],\label{id_25}\\ \
[a\cdot b,c]&=&a\cdot [b,c]+b\cdot [a,c]=[b,a\cdot c]+[a,b\cdot c],\label{id_26}
\eee
where $a\cdot b=\tfrac12(ab+ba)$, $ab$ is the associative product in $M_2(F)$.
\begin{lem}\label{lem4}
The following identities hold in $M_2(F)^{(+)}$ for $a,b,c\in H_2(F)$ 
\bes
a\cdot k&=&\tfrac12t(a)k,\\ \
[a,k]\cdot k&=&0,\\ \
[[k,a],k]&=&4a-2t(a),\\ \
[a,[k,b]]&=&(t(a)t(b)-2t(a\cdot b))k,\\ \
[k,a]\cdot [k,b]&=&2t(a\cdot b)-t(a)t(b),\\ \
[[k,a],[k,b]]&=&[[[k,a],k],b]=4[a,b],\\ \
[[a,[k,b]],k]&=&0,\\
([a,b],k,k)_+&=&0,\\
2a\cdot b&=&t(a)b+t(b)a-t(a)t(b)+t(a\cdot b),
\ees
where $t(a)\in F$ means the trace of the element $a$ of the quadratic algebra $H_2(F)$. 
\end{lem}
The proof follows from identity \eqref{id_25}, the multiplication table in $M_2(F)$, and the quadraticity of $H_2(F)$. 

The following theorem is a reminiscence of old results by P. Cohn \cite{C}.
\begin{thm}\label{thm2}
Assume that the bracket $\{x,y\}$ is  defined on the whole Jordan algebra $S=Z+N$. 
Then there exists an associative algebra with involution $(A,\star)$ such that  the algebra $J=H_2(F)\otimes Z+Fk\otimes N$ with the multiplication defined above is 
isomorphic to the Jordan algebra $H_2(A,\star)$ of $\star$-hermitian $2\times 2$ matrices over $(A,\star)$.
Conversly,  a Jordan algebra $H_2(A,*)$ of $2\times 2$ $*$-hermitian matrices over a unital associative algebra with involution $(A,*)$ has form  described in the theorem, with the bracket $\{,.\}$ defined on the whole $S$.
\end{thm}
\Proof
Define on the vector space $S$ a new multiplication
\bes
x*y=xy+2\{x,y\}
\ees
and prove that it is associative.  We have
\bes
(x*y)*z&=&(xy+2\{x,y\})*z=(xy)z+2\{xy,z\}+2\{x,y\}z+4\{\{x,y\},z\}\\
&=&(xy)z+2x\{y,z\}+2\{x,z\}y+2\{x,y\}z+4\{\{x,y\},z\},\\[1mm]
x*(y*z)&=&x*(yz+2\{y,z\})=x(yz)+2\{x,yz\}+2x\{y,z\}+4\{x,\{y,z\}\}\\
&=&x(yz)+2y\{x,z\}+2\{x,y\}z+2x\{y,z\}+4\{x,\{y,z\}\}.
\ees
Now, in view of \eqref{id_mult} and \eqref{id_cicl},
\bes
(x*y)*z-x*(y*z)=(x,y,z)-(x,z,y)-(y,x,z)=(x,y,z)+(y,z,x)+(z,x,y)=0.
\ees
Denote  the obtained associative algebra by $A$.  It is clear that $A^{(+)}=S$.  Consider the map $\star :A\rightarrow A,\ (z+n)^{\star}=z-n$  and prove that it is an involution in $A$.
We have for $z,u\in Z,\, n,m\in N$
\bes
((z+n)^{\star})^{\star}&=&(z-n)^{\star}=z+n,\\
((z+n)*(u+m))^{\star}&=&((z+n)(u+m))^{\star}+2\{(z+n),(u+m)\}^{\star}\\
&=&zu-un-zm+nm+2(\{n,u\}+\{z,m\}-\{z,u\}-\{n,m\})\\
&=&(u-m)(z-n)+2\{u-m,z-n\}=(u+m)^{\star}*(z+n)^{\star}.
\ees
Consider the associative algebra $M_2(A)\cong M_2(F)\otimes A$ and the adjoint Jordan algebra $M_2(A)^{(+)}$.
The product in $M_2(A)^{(+)}$ has form
\bes
(a\otimes x)\cdot(b\otimes y)&=&\tfrac12 ((a\otimes x)(b\otimes y)+(b\otimes y)(a\otimes x))\\
&=&\tfrac12(ab\otimes x*y+ba\otimes y*x)=\tfrac12(ab\otimes (xy+2\{x,y\})+ba\otimes (yx+2\{y,x\}))\\
&=&a\cdot b\otimes xy+[a,b]\otimes\{x,y\}.
\ees
Therefore, $J=H_2(F)\otimes Z+Fk\otimes N$ is a subalgebra of the algebra $M_2(A)^{(+)}$.  Moreover, let us show that $J$ coincides with the set of $\star$-hermitian matrices in $M_2(A)$.
In fact,  the involution   $^{-}$ induced by $\star$ in $M_2(A)=M_2(F)\otimes A$  has form 
\bes
\overline{x\otimes y}=x^t\otimes y^{\star},
\ees
where $x\rightarrow x^t$ is the matrix transposition. Therefore,  we have for any $a\in H_2, \, z\in Z, \,n\in N$
\bes
\overline{a\otimes z+k\otimes n}=a^t\otimes z^{\star}+k^t\otimes n^{\star}=a\otimes z+k\otimes n,
\ees
hence $J\subseteq H_2(A,\star)$.  Reciprocally, it is easy to see that any $\star$-hermitian matriz from $M_2(A)$ belongs to $J$.

\smallskip
Now, let $J=H_2(A,*)$ be the algebra of $2\times 2$ $*$-hermitian matrices over a unital associative algebra with involution $(A,*)$.  Decompose the vector space $A$ into the direct sum of $*$-symmetric and $*$-skewsymmetric elements  $A=H\oplus K$,  then $A^{(+)}=H+K$ is a $\Z_2$-graded Jordan algebra with $A^{(+)}_0=H$ and $A^{(+)}_1=K$.  For any $x,y\in A^{(+)}$,  put 
$\{x,y\}=\tfrac14(xy-yx)$, where $xy$ is the product in $A$.  Then $\{H,H\}\subseteq K,\, \{H,K\}\subseteq H,\, \{K,K\}\subseteq K$,  and identities \eqref{id_21'},\eqref{id_22} verify in $A^{(+)}$. 
Moreover, it is easy to check that $J$ is isomorphic to the algebra $H_2(F)\otimes H+Fk\otimes K$ with the product defined in the theorem. Hence we may take $S=A^{(+)}, \, Z=H,\, N=K$.
\ctd

\begin{rem}
Observe that the final part of the theorem is true for any $n$: the Jordan algebra $H_n(A,*)$ of $n\times n$ $*$-symmetric matrices over an associative algebra with involution $(A,*)$ can also be presented as $H_n(F)\otimes H+K_n(F)\otimes K$ with the product defined in the theorem, where $H_n(F)$ and $K_n(F)$ denote the spaces of symmetric and antisymmetric $n\times n$ matrices, respectively.
\end{rem}

\begin{Cor}\label{cor1}
Let $J=J_0+J_1$ be a $Z_2$-graded Jordan algebra with an  odd anticommutative bracket $\{.,.\}$.  
Consider the vector space  $A=H_2(F)\otimes J_0+Fk\otimes J_1$,  with the multiplication defined by 
\bee\label{id_mult1}
(a\otimes x)(b\otimes y)=a\cdot b\otimes xy+ [a,b]\otimes \{x,y\}.
\eee
Then, for any $a,b,c,d\in A$,  the expression $(ab,c,d)+(ad,c,b)+(bd,c,a)$ lies in the ideal $I$ of $A$ generated by the set
\bee\label{ideal}
 (x,y,z)-4\{y,\{x,z\}\}, \,\{xy,z\}-x\{y,z\}+y\{x,z\}, \,\{xy,z\}-\{y,xz\}+\{x,yz\},
\eee
where $x,y,z\in J$.
\end{Cor}

\section{The proof of the main theorem}

Let us finally prove

\begin{thm}\label{thm_2}
Let $J=J_0+J_1$ be a $Z_2$-graded Jordan algebra with a partial  odd anticommutative bracket $\{.,.\}$
\bes
\{.,.\}:J_i\times J_j\rightarrow J_{i+j+1},\   (i,j)\neq (1,1),
\ees
that satisfies identities 
\begin{equation*}
    \begin{split}
        (x,y,z)&=4\{y,\{x,z\}\}, \text{ for any $x, y, z \in J$,}\\
  \{xz,y\}&=x\{z,y\}+z\{x,y\}=\{z,xy\}+\{x,zy\}, \text{ for any $x, y, z \in J$,}\\
  (\{x,y\},&z,t)+(\{t,x\},z,y)+(\{y,t\},z,x)=0 \text{ for any $x, y \in J_0$ and $z, t \in J_1$},
    \end{split}
\end{equation*}
 and
always when all involved products are defined.
Then the vector space  $A=H_2(F)\otimes J_0+Fk\otimes J_1$,  with the multiplication defined by 
\eqref{id_mult1},  is a Jordan algebra.
\end{thm}
\Proof 
We have to prove that for any $a\otimes z,\, b\otimes u,\, c\otimes v,\, d\otimes w\in A$ the Jordan identity holds 
\begin{equation}\label{id_x}
\begin{split}
0&=(a\otimes z\cdot b\otimes u,c\otimes v,d\otimes w)+(a\otimes z\cdot d\otimes w,c\otimes v,b\otimes u)\\
&+(b\otimes u\cdot d\otimes w,c\otimes v,a\otimes z).
\end{split}
\end{equation}

After calculating the products in \eqref{id_x},  the right side of it may be written as $S=S_1+S_2+S_3$, where
\bes
&S_1=((a\cdot b)\cdot c)\cdot d\otimes (zu\cdot v)w+[(a\cdot b)\cdot c,d]\otimes \{zu\cdot v,w\}+&\\
&+[a\cdot b,c]\cdot d\otimes \{zu,v\}w+[[a\cdot b,c],d]\otimes\{\{zu,v\},w\}+&\\
&+([a,b]\cdot c)\cdot d\otimes (\{z,u\}v)w+[[a,b]\cdot c,d]\otimes \{\{z,u\}v,w\}+&\\
&+[[a,b],c]\cdot d\otimes \{\{z,u\},v\}w+[[[a,b],c],d]\otimes\{\{\{z,u\},v\},w\}-&\\
&-(a\cdot b)\cdot (c\cdot d)\otimes zu\cdot vw-[a\cdot b,c\cdot d]\otimes \{zu,vw\}-&\\
&-(a\cdot b)\cdot [c,d]\otimes (zu)\{v,w\}-[a\cdot b,[c,d]]\otimes \{zu,\{v,w\}\}-&\\
&-[a,b]\cdot (c\cdot d)\otimes \{z,u\}(vw)-[[a,b],c\cdot d]\otimes \{\{z,u\},vw\}-&\\
&-[a,b]\cdot [c,d]\otimes \{z,u\}\{v,w\}-[[a,b],[c,d]]\otimes \{\{z,u\},\{v,w\}\},&
\ees
\bes
&S_2=((a\cdot d)\cdot c)\cdot b\otimes (zw\cdot v)u+[(a\cdot d)\cdot c,b]\otimes \{zw\cdot v,u\}+&\\
&+[a\cdot d,c]\cdot b\otimes \{zw,v\}u+[[a\cdot d,c],b]\otimes\{\{zw,v\},u\}+&\\
&+([a,d]\cdot c)\cdot b\otimes (\{z,w\}v)u+[[a,d]\cdot c,b]\otimes \{\{z,w\}v,u\}+&\\
&+[[a,d],c]\cdot b\otimes \{\{z,w\},v\}u+[[[a,d],c],b]\otimes\{\{\{z,w\},v\},u\}-&\\
&-(a\cdot d)\cdot (c\cdot b)\otimes zw\cdot vu-[a\cdot d,c\cdot b]\otimes \{zw,vu\}-&\\
&-(a\cdot d)\cdot [c,b]\otimes (zw)\{v,u\}-[a\cdot d,[c,b]]\otimes \{zw,\{v,u\}\}-&\\
&-[a,d]\cdot (c\cdot b)\otimes \{z,w\}(vu)-[[a,d],c\cdot b]\otimes \{\{z,w\},vu\}-&\\
&-[a,d]\cdot [c,b]\otimes \{z,w\}\{v,u\}-[[a,d],[c,b]]\otimes \{\{z,w\},\{v,u\}\},&
\ees
\bes 
&S_3=((b\cdot d)\cdot c)\cdot a\otimes (uw\cdot v)z+[(b\cdot d)\cdot c,a]\otimes \{uw\cdot v,z\}+&\\
&+[b\cdot d,c]\cdot a\otimes \{uw,v\}z+[[b\cdot d,c],a]\otimes\{\{uw,v\},z\}+&\\
&+([b,d]\cdot c)\cdot a\otimes (\{u,w\}v)z+[[b,d]\cdot c,a]\otimes \{\{u,w\}v,z\}+&\\
&+[[b,d],c]\cdot a\otimes \{\{u,w\},v\}z+[[[b,d],c],a]\otimes\{\{\{u,w\},v\},z\}-&\\
&-(b\cdot d)\cdot (c\cdot a)\otimes uw\cdot vz-[b\cdot d,c\cdot a]\otimes \{uw,vz\}-&\\
&-(b\cdot d)\cdot [c,a]\otimes (uw)\{v,z\}-[b\cdot d,[c,a]]\otimes \{uw,\{v,z\}\}-&\\
&-[b,d]\cdot (c\cdot a)\otimes \{u,w\}(vz)-[[b,d],c\cdot a]\otimes \{\{u,w\},vz\}-&\\
&-[b,d]\cdot [c,a]\otimes \{u,w\}\{v,z\}-[[b,d],[c,a]]\otimes \{\{u,w\},\{v,z\}\}.&
\ees

The algebra  $A$ also has a $Z_2$-grading, and we may assume that the arguments in \eqref{id_x} are homogeneous.
Moreover,  without loss of generality,  we may assume that $J$  is a free algebra with the set of even generators $X$ and odd generators $Y$, and that the elements $z,u,v,w\in X\cup Y$.

\smallskip
Consider  first the case when $|\{z,u,v,w\}\cap Y|\leq 1$.  By Corollary \ref{cor1},  $S=\sum u_iW_i$, where $u_i$ are elements of $J$ of type \eqref{ideal} and $W_i$ some multiplication operators in $A$.  By homogenity,  among the arguments of every $u_i$ in this expression there is at most one odd element.  Therefore,   $u_i=0$  by the conditions of the theorem, and hence $S=0$.

\smallskip
Let  now $|\{z,u,v,w\}\cap Y|=2.$ We have the two cases.  First,  let $z,u\in Y,\, v,w\in X$.   In this case, due to the symmetry (and further linearization on $z$,)  it suffices to calculate $S$ for $u=z,\ a=b=k$.  We have by Lemma \ref{lem4} and identities \eqref{id_21'},\eqref{id_22}

\bes
S_1&=(k^2\cdot c)\cdot d\otimes (z^2\cdot v)w+[k^2\cdot c,d]\otimes \{z^2v,w\}+&\\
&+[k^2,c]\cdot d\otimes \{z^2,v\}w+[[k^2,c],d]\otimes\{\{z^2,v\},w\}+&\\
&-k^2\cdot (c\cdot d)\otimes z^2\cdot vw-[k^2,c\cdot d]\otimes \{z^2,vw\}-&\\
&-(k^2)\cdot [c,d]\otimes (z^2)\{v,w\}-[k^2,[c,d]]\otimes \{z^2,\{v,w\}\}&\\
&=-(c\cdot d)\otimes (z^2\cdot v)w-[c,d]\otimes \{z^2v,w\}&\\
&+(c\cdot d)\otimes z^2\cdot vw+[c,d]\otimes (z^2)\{v,w\}&\\
&=-(c\cdot d)\otimes (z^2, v,w)-[c,d]\otimes v\{z^2,w\},
\ees

\bes
&S_2=S_3=((k\cdot d)\cdot c)\cdot k\otimes (zw\cdot v)z+[(k\cdot d)\cdot c,k]\otimes \{zw\cdot v,z\}+&\\
&+[k\cdot d,c]\cdot k\otimes \{zw,v\}z+[[k\cdot d,c],k]\otimes\{\{zw,v\},z\}+&\\
&+([k,d]\cdot c)\cdot k\otimes (\{z,w\}v)z+[[k,d]\cdot c,k]\otimes \{\{z,w\}v,z\}+&\\
&+[[k,d],c]\cdot k\otimes \{\{z,w\},v\}z+[[[k,d],c],k]\otimes\{\{\{z,w\},v\},z\}-&\\
&-(k\cdot d)\cdot (c\cdot k)\otimes zw\cdot vz-[k\cdot d,c\cdot k]\otimes \{zw,vz\}-&\\
&-(k\cdot d)\cdot [c,k]\otimes (zw)\{v,z\}-[k\cdot d,[c,k]]\otimes \{zw,\{v,z\}\}-&\\
&-[k,d]\cdot (c\cdot k)\otimes \{z,w\}(vz)- [[k,d],c\cdot k]\otimes \{\{z,w\},vz\}-&\\
&-[k,d]\cdot [c,k]\otimes \{z,w\}\{v,z\}-[[k,d],[c,k]]\otimes \{\{z,w\},\{v,z\}\}.&\\
\ees
By Lemma \ref{lem4} we have
\bes
&((k\cdot d)\cdot c)\cdot k=-\tfrac14 t(d)t(c)&\\
&[(k\cdot d)\cdot c,k]=\tfrac14 t(d)t(c)[k,k]=0&\\
&[k\cdot d,c]\cdot k=\tfrac12 t(d)[k,c]\cdot k=0&\\
&[[k\cdot d,c],k]=\tfrac12 t(d)[[k,c],k]&\\
&([k,d]\cdot c)\cdot k=[k\cdot c,d]\cdot k-([c,d]\cdot k)\cdot k=&\\
&\tfrac12 t(c)[k,d]\cdot k-[c,d]k^2=[c,d]&\\
&[[k,d]\cdot c,k]=[[k\cdot c,d],k]-[[c,d]\cdot k,k]=\tfrac12 t(c)[[k,d],k]&\\
&[[k,d],c]\cdot k=[[k,d]\cdot k,c]-[k,d]\cdot [k,c]=-[k,c]\cdot [k,d]&\\
&[[[k,d],c],k]=[-(t(c)t(d)+2t(c\cdot d))k,k]=0&\\
&-(k\cdot d)\cdot (c\cdot k)=\tfrac14 t(d)t(c)&\\
&-[k\cdot d,c\cdot k]= -\tfrac14 t(d)t(c)[k,k]=0&\\
&-(k\cdot d)\cdot [c,k]=-\tfrac12 t(d)k\cdot [c,k]=0&\\
&-[k\cdot d,[c,k]]=-\tfrac12 t(d)[[k,c],k]&\\
&-[k,d]\cdot (c\cdot k)=-\tfrac12 t(c)[k,d]\cdot k=0&\\
&-[[k,d],c\cdot k]=-\tfrac12 t(c)[[k,d],k]&\\
&-[[k,d],[c,k]]=4[d,c].&
\ees
Hence
\bes
&S_2=S_3=-\tfrac14 t(d)t(c)\otimes (zw\cdot v)z+\tfrac12 t(d)[[k,c],k]\otimes\{\{zw,v\},z\}+&\\
&+[c,d]\otimes (\{z,w\}v)z+\tfrac12 t(c)[[k,d],k]\otimes \{\{z,w\}v,z\}-&\\
&-[k,c]\cdot [k,d]\otimes \{\{z,w\},v\}z-\tfrac14 t(d)t(c)\otimes zw\cdot vz-&\\
&-\tfrac12 t(d)[[k,c],k]\otimes \{zw,\{v,z\}\}-\tfrac12 t(c)[[k,d],k]\otimes \{\{z,w\},vz\}-&\\
&-[k,d]\cdot [c,k]\otimes \{z,w\}\{v,z\}+4[d,c]\otimes \{\{z,w\},\{v,z\}\}&\\
&=-\tfrac14 t(d)t(c)\otimes (zw,v,z)+\tfrac12 t(d)[[k,c],k]\otimes(\{\{zw,v\},z\}-\{zw,\{v,z\}\})&\\
&+\tfrac12 t(c)[[k,d],k]\otimes (\{\{z,w\}v,z\}-\{\{z,w\},vz\})&\\
&+[c,d]\otimes ((\{z,w\}v)z-4\{\{z,w\},\{v,z\}\})&\\
&-[k,c]\cdot [k,d]\otimes (\{\{z,w\},v\}z-\{z,w\}\{v,z\}).
\ees
We have
\bes
&-\tfrac14  (zw,v,z)=-\{v,\{zw,z\}\}=-\tfrac12\{v,\{w,z^2\}\},&\\
&\tfrac12\{\{zw,v\},z\}-\{zw,\{v,z\}\}=\tfrac18 (-(zw,z,v)-(v,zw,z))=&\\
&\stackrel{\eqref{id_cicl}}=\tfrac18 (z,v,zw)=\tfrac12\{v,\{z,zw\}\}=\tfrac14\{v,\{z^2,w\}\},&\\
&\tfrac12\{\{z,w\}v,z\}-\{\{z,w\},vz\}\stackrel{\eqref{id_22}}=\tfrac12\{v,z\{z,w\}\}=\tfrac14\{v,\{z^2,w\}\},&\\
&(\{z,w\}v)z-4\{\{z,w\},\{v,z\}\}=(\{z,w\}v)z-(v,\{z,w\},z)=v(\{z,w\}z)=\tfrac12 v\{z^2,w\},&\\
&\{\{z,w\},v\}z-\{z,w\}\{v,z\}=\{z\{z,w\},v\}=\tfrac12 \{\{z^2,w\},v\}.&
\ees
Therefore,
\bes
&S_2=S_3=\tfrac14(-2t(c)t(d)-t(d)[[k,c],k]- t(c)[[k,d],k]-2[k,c]\cdot [k,d])\otimes\{\{z^2,w\},v\})&\\
&+\tfrac12[c,d]\otimes v\{z^2,w\}.&
\ees
By Lemma \ref{lem4},
\bes
&-2t(c)t(d)-t(d)[[k,c],k]- t(c)[[k,d],k]-2[k,c]\cdot [k,d]=&\\
&-2t(c)t(d)-t(d)(4c-2t(c))-t(c)(4d-2t(d))-4t(c\cdot d)+2t(c)t(d)=&\\
&=4(t(c)t(d)-t(c)d-t(d)c)-t(c\cdot d))=&\\
&=-8c\cdot d.
\ees
Hence 
\bes
S_2&=&-2c\cdot d\otimes \{\{z^2,w\},v\}+\tfrac12[c,d]\otimes v\{z^2,w\}\\
&\stackrel{\eqref{id_21'}}=&\tfrac12(c\cdot d\otimes (z^2,v,w)+[c,d]\otimes v\{z^2,w\}),
\ees
and $S_1+2S_2=0$. 

 In the second case we have $z,v\in Y$, hence we can put $a=c=k$. Then
\bes 
&S_1=((k\cdot b)\cdot k)\cdot d\otimes (zu\cdot v)w+[[k,b],k]\cdot d\otimes \{\{z,u\},v\}w&\\
&+[[[k,b],k],d]\otimes\{\{\{z,u\},v\},w\}-(k\cdot b)\cdot (k\cdot d)\otimes zu\cdot vw&\\
&-[k\cdot b,[k,d]]\otimes \{zu,\{v,w\}\}-[[k,b],k\cdot d]\otimes \{\{z,u\},vw\}-&\\
&-[k,b]\cdot [k,d]\otimes \{z,u\}\{v,w\}-[[k,b],[k,d]]\otimes \{\{z,u\},\{v,w\}\},&
\ees
\bes
&S_2=((k\cdot d)\cdot k)\cdot b\otimes (zw\cdot v)u+[[k,d],k]\cdot b\otimes \{\{z,w\},v\}u&\\
&+[[[k,d],k],b]\otimes\{\{\{z,w\},v\},u\}-(k\cdot d)\cdot (k\cdot b)\otimes zw\cdot vu&\\
&-[k\cdot d,[k,b]]\otimes \{zw,\{v,u\}\}-[[k,d],k\cdot b]\otimes \{\{z,w\},vu\}&\\
&-[k,d]\cdot [k,b]\otimes \{z,w\}\{v,u\}-[[k,d],[k,b]]\otimes \{\{z,w\},\{v,u\}\},&
\ees
\bes
&S_3=((b\cdot d)\cdot k)\cdot k\otimes (uw\cdot v)z+[[b\cdot d,k],k]\otimes\{\{uw,v\},z\}&\\
&+([b,d]\cdot k)\cdot k\otimes (\{u,w\}v)z&\\
&-(b\cdot d)\cdot (k\cdot k)\otimes uw\cdot vz-[b,d]\cdot (k\cdot k)\otimes \{u,w\}(vz).&
\ees
We have
\bes
&((k\cdot b)\cdot k)\cdot d=-\tfrac12 t(b)d,&\\
&[[k,b],k]\cdot d=(4b-2t(b))\cdot d =4b\cdot d-2t(b)d=&\\
&=2(t(d)b-t(b)t(d)+t(b\cdot d)),&\\
&[[[k,b],k],d]=[4b-2t(b),d]=4[b,d],&\\
&-(k\cdot b)\cdot (k\cdot d)=\tfrac14 t(b)t(d),&\\
&-[k\cdot b,[k,d]]=\tfrac12 t(b)[[k,d],k]=t(b)(2d-t(d))=2t(b)d-t(b)t(d),&\\
&-[[k,b],k\cdot d]=-\tfrac12 t(d)[[k,b],k]=-2t(d)b+t(b)t(d),&\\
&-[k,b]\cdot [k,d]=-2t(b\cdot d)+t(b)t(d),&\\
&-[[k,b],[k,d]]=-4[b,d].&
\ees
Note that in $S_2$ the first terms of the tensor products of each summand are the same as in $S_1$ by changinng the roles of $b$ and $d$. The same, in $S_2$, the second terms in the tensor products are like those of $S_1$ interchanging $u$ and $w$.
\bes
&((b\cdot d)\cdot k)\cdot k=-\tfrac12 t(b\cdot d),&\\
&[[b\cdot d,k],k]=\tfrac12(t(b)[[d,k],k]+t(d)[[b,k],k])=&\\
&=-2(t(b)d+t(d)b-t(b)t(d)),&\\
&([b,d]\cdot k)\cdot k=-[b,d],&\\
&-(b\cdot d)\cdot (k\cdot k)=b\cdot d=\tfrac12 (t(b)d+t(d)b-t(b)t(d)+t(b\cdot d)),&\\
&-[b,d]\cdot (k\cdot k)=[b,d].&
\ees
Hence we may  write
\bes
&S_1+S_2+S_3=&\\
&t(b)d\otimes A_1+t(d)b\otimes  A_2+t(b)t(d)\otimes A_3+t(b\cdot d)\otimes A_4+[b,d]\otimes A_5.&
\ees
We have
\bes
&A_1=-\tfrac12 (zu\cdot v)w+2\{zu,\{v,w\}\}+2 \{\{z,w\},v\}u-2 \{\{z,w\},vu\}&\\
&-2\{\{uw,v\},z\}+\tfrac12 uw\cdot vz\stackrel{\eqref{id_21'},\eqref{id_4}}=&\\
&=\tfrac12 ((v,zu,w)+v(z,u,w)-(v,z,uw)+(v,z,u)w-(vz,u,w))\stackrel{\eqref{Teich}}=0,&
\ees
Similarly,  $A_2=0$.  Furthermore,
\bes
&A_3=-2\{\{z,u\},v\}w-2 \{\{z,w\},v\}u+\tfrac14 zu\cdot vw+\tfrac14 zw\cdot vu- \{zu,\{v,w\}\}&\\
& - \{zw,\{v,u\}\}+\{\{z,u\},vw\}+ \{\{z,w\},vu\}+2\{\{uw,v\},z\}&\\
&-\tfrac12 uw\cdot vz+\{z,u\}\{v,w\}+\{z,w\}\{v,u\}&\\
&\stackrel{\eqref{id_21'},\eqref{id_22}}=\tfrac14((zu\cdot vw+zw\cdot vu-2uw\cdot vz)+2(z,v,u)w+2(z,v,w)u&\\
&-(v,zu,w)-(v,zw,u)-(z,vw,u)-(z,vu,w)-2(uw,z,v))+&\\
&+\{z,u\{v,w\}\}+\{z,w\{v,u\}\}-\{z,\{v,w\}\}u-\{z,\{v,u\}\}w&\\
&=\tfrac14(-(zu,w,v)-(zw,u,v)+2(uw,z,v)+(z,u,w)v+(z,w,u)v&\\
&+2(z,v,u)w+2(z,v,w)u-(v,zu,w)-&\\
&-(v,zw,u)-(z,vw,u)-(z,vu,w)-2(uw,z,v))+&\\
&+\{z,\{v,uw\}\}-\{z,\{v,w\}\}u-\{z,\{v,u\}\}w=&\\
&\stackrel{\eqref{id_cicl},\eqref{id_2},\eqref{id_4},\eqref{id_21'}}=\tfrac14((w,v,zu)+(u,v,zw)+(z,v,u)w+(z,v,w)u+&\\
&+(v,z,uw)-(v,z,w)u-(v,z,u)w)=&\\
&\stackrel{\eqref{id_3}}=\tfrac14((uw,v,z)+(v,z,uw)+(z,u,v)w+(z,w,v)u)=&\\
&\stackrel{\eqref{id_cicl},\eqref{id_4}}=\tfrac14(-(z,uw,v)+(z,u,v)w+(z,w,v)u)=0.&
\ees
Now,
\bes
&A_4= 2\{\{z,u\},v\}w-2 \{z,u\}\{v,w\}+2\{\{z,w\},v\}u-2\{z,w\}\{v,u\}&\\
&-\tfrac12 (uw\cdot v)z +\tfrac12 uw\cdot vz=&\\
&=2\{\{z,u\}w,v\}+2\{\{z,w\}u,v\}-\tfrac12 (uw,v,z)=&\\
&=2\{\{z,uw\},v\}-\tfrac12 (uw,v,z)=&\\
&=\tfrac12( -(z,v,uw) -(uw,v,z))=0.
\ees
Finally,
\bes 
&A_5=4(\{\{\{z,u\},v\},w\}-\{\{z,u\},\{v,w\}\}-\{\{\{z,w\},v\},u\}&\\
&+\{\{z,w\},\{v,u\}\})-(\{u,w\}v)z+\{u,w\}(vz)=&\\
&=-(\{z,u\},w,v)-(v,\{z,u\},w)+(\{z,w\},u,v)+(v,\{z,w\},u)-(\{u,w\},v,z)=&\\
&=(w,v,\{z,u\})-(u,v,\{z,w\})+(z,v,\{u,w\})\stackrel{\eqref{id_22'}}=0.&
\ees
\smallskip

Let now  $|\{z,u,v,w\}\cap Y|=3$.  Again, we have the two cases.  Assume first that $z,u,w\in Y$. Then we may put $a=b=d=k$,  hence
\bes
&S_1=((k\cdot k)\cdot c)\cdot k\otimes (zu\cdot v)w+[(k\cdot k)\cdot c,k]\otimes \{zu\cdot v,w\}+&\\
&-(k\cdot k)\cdot (c\cdot k)\otimes zu\cdot vw-&\\
&-(k\cdot k)\cdot [c,k]\otimes (zu)\{v,w\}&\\
&=-\tfrac12 t(c)k\otimes (zu\cdot v)w-[c,k]\otimes \{zu\cdot v,w\}&\\
&+\tfrac12 t(c)k\otimes zu\cdot vw+[c,k]\otimes (zu)\{v,w\}&\\
&=-\tfrac12 t(c)k\otimes (zu,v,w)-[c,k]\otimes v\{zu,w\}.
\ees
Similarly, 
\bes
S_2&=&-\tfrac12 t(c)k\otimes (zw,v,u)-[c,k]\otimes v\{zw,u\},\\
S_3&=&-\tfrac12 t(c)k\otimes (uw,v,z)-[c,k]\otimes v\{uw,z\},\\
\ees
and
\bes
S&=&-\tfrac12 t(c)k\otimes ((zu,v,w)+(zw,v,u)+(uw,v,z))\\ \ 
&-&[c,k]\otimes v( \{zu,w\}+\{zw,u\}+\{uw,z\} )=0.
\ees
In the second case we may assume that $z,u,v\in Y$.  Hence we may assume that $a=b=c=k$, and then
\bes
&S_1=((k\cdot k)\cdot k)\cdot d\otimes (zu\cdot v)w+[(k\cdot k)\cdot k,d]\otimes \{zu\cdot v,w\}&\\
&-(k\cdot k)\cdot (k\cdot d)\otimes zu\cdot vw-(k\cdot k)\cdot [k,d]\otimes (zu)\{v,w\}&\\
&=-\tfrac12 t(d)k\otimes (zu\cdot v)w-[k,d]\otimes \{zu\cdot v,w\}&\\
&\tfrac12 t(d)k\otimes zu\cdot vw+[k,d]\otimes (zu)\{v,w\}&\\
&=-\tfrac12 t(d)k\otimes (zu, v,w)-[k,d]\otimes v\{zu,w\},&
\ees
\bes
&S_2=((k\cdot d)\cdot k)\cdot k\otimes (zw\cdot v)u+[[[k,d],k],k]\otimes\{\{\{z,w\},v\},u\}&\\
&-(k\cdot d)\cdot (k\cdot k)\otimes zw\cdot vu-[k,d]\cdot (k\cdot k)\otimes \{z,w\}(vu)&\\
&=-\tfrac12 t(d)k\otimes (zw\cdot v)u+4[d,k]\otimes\{\{\{z,w\},v\},u\}&\\
&+\tfrac12 t(d)k\otimes zw\cdot vu+[k,d]\otimes \{z,w\}(vu)&\\
&=-\tfrac12 t(d)k\otimes (zw,v,u)+[d,k]\otimes (4\{\{\{z,w\},v\},u\}- \{z,w\}(vu)),&
\ees
\bes 
&S_3=((k\cdot d)\cdot k)\cdot k\otimes (uw\cdot v)z+[[[k,d],k],k]\otimes\{\{\{u,w\},v\},z\}&\\
&-(k\cdot d)\cdot (k\cdot k)\otimes uw\cdot vz-[k,d]\cdot (k\cdot k)\otimes \{u,w\}(vz)&\\
&=-\tfrac12 t(d)k\otimes (uw\cdot v)z+4[d,k]\otimes\{\{\{u,w\},v\},z\}&\\
&+\tfrac12 t(d)k\otimes uw\cdot vz+[k,d]\otimes \{u,w\}(vz)&\\
&=-\tfrac12 t(d)k\otimes (uw,v,z)+[d,k]\otimes(4\{\{\{u,w\},v\},z\}- \{u,w\}(vz)).&
\ees
Therefore,
\bes 
S&=&-\tfrac12 t(d)k\otimes ((zu, v,w)+(zw,v,u)+(uw,v,z))\\
&+&[d,k]\otimes (v\{zu,w\}+4\{\{\{z,w\},v\},u\}- \{z,w\}(vu)+4\{\{\{u,w\},v\},z\}- \{u,w\}(vz))\\
&=&[d,k]\otimes (v(z\{u,w\})+v(u\{z,w\})- \{z,w\}(vu)- \{u,w\}(vz))\\
&+&4\{\{\{z,w\},v\},u\}+4\{\{\{u,w\},v\},z\})\\
&=&[d,k]\otimes( (\{u,w\},z,v)+(\{z,w\},u,v)-(\{z,w\},u,v)-(\{u,w\},z,v))=0.
\ees
\smallskip

Finally, let $z,u,v,w\in Y$. Then we may take $a=b=c=d=k$, hence
\bes
&S_1=((k\cdot k)\cdot k)\cdot k\otimes (zu\cdot v)w&\\
&-(k\cdot k)\cdot (k\cdot k)\otimes zu\cdot vw=1\otimes (zu,v,w).&
\ees
Similarly, $S_2=1\otimes (zw,v,u),\ S_3=1\otimes (uw,v,z)$, and eventually
\bes
S=1\otimes ((zu,v,w)+(zw,v,u)+(uw,v,z))=0.
\ees
This proves the theorem.

\ctd

The bracket $\{,\}$ may really not be defined on the whole algebra $S$.  The most evident such case is given by the following
\begin{prop}\label{prop_5}
The bracket $\{,\}$ is trivial, that is, equal to zero on $S=Z+N$, (and not defined on $N$) if and only if the algebra $J$ is isomorfic to an algebra $J=A+V$ of a bilinear form  $f:V\otimes V\rightarrow A$ on an associative commutative module $V$ over a unital commutative associative algebra $A$ such that there exist $u,v\in V$ with $f(u,u)=f(v,v)=1,\, f(u,v)=0$.
\end{prop}
\Proof
Let first the bracket be trivial on $S$.  Then the algebra $Z$ is associative,  since $(Z,Z,Z)=\{Z,\{Z,Z\}\}=0$. Moreover,  the $Z$-bimodule $N$ is associative as welll since we have $(N,Z,Z)=(Z,Z,N)=\{Z,\{N,Z\}\}=0$ and by \eqref{id_cicl} $(Z,N,Z)\subseteq (Z,Z,N)+(N,Z,Z)=0$.  Now, put $A=Z, \,V=N+Z(1-2e)+Zh, \,f(u,v)=uv$, then $S\cong A+V$ as an algebra of bilinear form, and $f(1-2e,h)=(1-2e)h=0, \, f(h,h)=f(1-2e,1-2e)=1.$ 

Conversely, let $J=A+V$ be the algebra of $A$-bilinear form $f$ on an $A$-bimodule $V$ and $u,v\in V$ be the elements satisfying the condition of the proposition.  Then the $F$-subspace $F1+Fu+Fv$ is a subalgebra of $J$ isomorfic to  $H_2(F)$,  and we have $Z=A$ and $N$ equal to the orthogonal complement of $Fv+Fu$ in $V$ with respect to $f$.

\ctd

It remains an open question whether there exist $H_2$-algebras different from algebras of type $H_2(A,*)$ and algebras of bilinear form.
Observe that such examples neither appear  in the classification of Osborn \cite{O} .

\smallskip

It follows from Zelmanov's classification theorem that every simple $H_2$-algebra is of type $H_2(A,*)$ or an algebra of bilinear form.
It would be interesting obtain a direct proof of this  fact.
It follows from our results that this fact is equivalent to the condition that for a simple $H_2$-algebra $J$ the corresponding  anticommutative odd bracket is either completely defined or trivial.

Another open question is the following: Is it true that any $H_2(F)$-algebra is special?.

\end{document}